\magnification=1200
\def\R{\hbox{$I\kern-3.5pt R$}}
\def\N{\hbox{$I\kern-3.5pt N$}}
\def\E{\hbox{$I\kern-3.5pt E$}}
\font\fp=cmr8

\def\b{\bigskip}
\def\m{\medskip}
\def\d{\displaystyle }

\def\P{\hbox{$I\kern-3.5pt P$}}
\def\Rm#1{\vskip 6 pt \noindent {\bf Remark #1.\ }}
\def\Pf{{\bf Proof: }}
\def\fin{ \hfill /// \vskip 12pt}

\def\Rm{\vskip 6 pt \noindent {\bf Remark.\ }}

\def\sk{\sum_{i=1}^k}
\def\s{\sum_{i=1}^n}
\hsize=15.5true cm
\vsize=23true cm
\def\pconv{\hbox{$p$-{\rm conv}\ }}
\def\qconv{\hbox{$q$-{\rm conv}\ }}
\def\conv{\hbox{{\rm conv}\ }}
\def \lx{\s\lambda_i x_i}
\def \my{\sum_{i=1}^m\mu_i y_i}
\def \l{\s\lambda_i^p}
\hfuzz=1pt

\null

Preliminary version.

\bigskip
\centerline {\bf The theorems of Caratheodory and Gluskin for $0<p<1$}

\vskip 6 pt
\centerline {by}
\vskip 6pt
\centerline { Jes\'us Bastero \footnote 
*{Partially supported  by Grant DGICYT PS 90-0120} 
, Julio Bernu\'es * and Ana Pe\~na \footnote
{**}{Supported by Grant DGA (Spain)}}

\vskip
6pt \centerline{\fp Departamento de Matem\'aticas. Facultad de Ciencias}
\centerline{\fp Universidad de Zaragoza}
\centerline{\fp 50009-Zaragoza (Spain)} 
\vskip 20pt

\beginsection 1. Introduction and notation.\par
\m

Throughout the paper $X$ will denote a real vector space and $p$ will be a real
number, $0<p<1$. A set $A\subseteq X$ is called $p$-convex if $\lambda x+\mu
y\in A$, whenever $ x,y\in A$,  and $ \lambda, \mu\geq 0 $, with
$\lambda^p+\mu^p=1$.  Given $A\subseteq X$, the $p$-convex hull of $A$ is
defined as the intersection of all $p$-convex sets that contain $A$. Such set is
denoted by $\pconv (A)$.

 A $p$-norm on $X$ is a map $\Vert\cdot\Vert\colon X\to\R$
verifying:
\item{(i)} $\Vert x\Vert\geq 0, \forall\, x\in X$ and $\Vert
x\Vert=0\Leftrightarrow x=0$.
\item{(ii)} $\Vert ax\Vert=|a|\,\Vert x\Vert, \forall\,a\in\R, x\in X$.
\item{(iii)} $\Vert x+y\Vert^p \leq \Vert x\Vert^p+\Vert y\Vert^p, \forall\,
x,y\in X$.

We will say that $(X,\Vert\cdot\Vert)$ is a $p$-normed space. The unit ball of
a $p$-normed space is a $p$-convex set and will be denoted by $B_X$.

We denote by ${\cal M}_n^p$ the class of all $n$-dimensional $p$-normed
spaces. If $X,Y \in {\cal M}_n^p $ the Banach-Mazur distance $d(X,Y)$ is
the infimun of the products $\Vert T \Vert\cdot \Vert T^{-1} \Vert$, where the
infimun is taken over all the isomorphisms $T$ from $X$ onto $Y$.

We shall use the notation and terminology commonly used in Banach space theory
as it appears in [Tmcz].

In this note we investigate some aspects of the local structure of finite
dimensional $p$-Banach spaces. The well known theorem of Gluskin gives a
sharp lower bound of the diameter of the Minkowski compactum. In [Gl] it is proved
that diam$({\cal M}_n^1)\geq cn$ for some absolute constant $c$. 
Our purpose is to study this problem in the $p$-convex setting. In [Pe], T. Peck
 gave an upper bound of the diameter of ${\cal M}_n^p$ namely, 
diam$({\cal M}_n^p)\leq n^{2/p-1}$. We will show that such bound is optimum. 

The method used by Gluskin to prove his result can be directed applied, with some
minor variations, to our case. At some point of the proof it is necessary to
find some  volumetric estimates for convex envelopes. In particular if $\{
P_i \}_{i=1}^m $ are $m$ points in the euclidean sphere in $ \R^n $  we need to
estimate from above $\d \left( {\vert \conv\{\pm P_i\} \vert \over \vert
B_{\ell_2^n} \vert} \right)^{1/n}$. Szarek, [Sz], and other authors gave the
 estimate $\d \left( {\vert \conv\{\pm P_i\} \vert
\over \vert B_{\ell_2^n} \vert} \right)^{1/n}\leq C mn^{-3/2}$ ($C$ is an
absolute constant).  Caratheodory's convexity theorem turned out to be an
important 
 ingredient of the proof. For $p<1$ we will proceed in this fashion
and so we will need to have the corresponding version of Caratheodory convexity
theorem.

 The main results of the paper are the following:\b

\proclaim Theorem 1. Let $ A\subseteq \R^n$ and $0<p<1$. For every  $x\in \pconv
({A}), x\not=0$ there exist linearly independent vectors $\{ P_1 \ldots
P_k\}\subseteq A$  with $k\le n$, such that $x\in \pconv{\{ P_1 \ldots P_k\}}$. 
Moreover, if  $0\in \pconv({A}) $, there exits $\{ P_1 \ldots P_k\} \subseteq A$
with $k\le n+1 $ such that $0\in \pconv{\{ P_1 \ldots P_k\}}$.

\m

\proclaim Theorem 2. Let  $0<p<1$. There exits a constant $C_p >0$ such that
for every $n\in \N $ 
$$   C_p n^{2/p-1}\leq {\rm diam}({\cal M}_n^p)\leq n^{2/p-1}.$$

\m

The first result can be viewed as the $p$-convex analogue of Caratheodory's
theorem. Apparently, the result for $p<1$ is stronger than the Caratheodory's one
in the sense that we get  $k\le n$ and only $k\le n+1 $ can be assured for
$p=1$ (see [Eg], pg 35). Proposition 3 ii) will show that this is not such since
vector 0 plays a particularly special role.

The second result is the analogue of Gluskin's theorem in the $p$-convex
setting, that is the diameter of the Minkowski compactum grows asimpthotically
like $n^{2/p-1}$.
\vskip 12pt

\noindent{\bf 2. Caratheodory's theorem for $p$-convex hulls.}
\vskip 12pt

 In this section we want to prove Theorem 1. We begin by recalling the
first properties of $p$-convex hulls. They are probably known but since we have
not found them in any reference we sketch their proofs.

\proclaim Proposition 3. For every $\emptyset\ne A\subseteq X$.
\item{i)} \ \ $\d \pconv (A)=
\left\{ \lx\mid \lambda_i\geq 0, \l=1,
x_i\in A, n\in \N\right\}$.
\medskip
\item{ii)}\vskip-15pt
$\eqalign{\d \,\,\pconv (A\cup\{ 0\}) &= \{ 0\}\bigcup \pconv
(A) \cr &= \left\{ \lx\mid \lambda_i\geq 0,
\l\leq1, x_i\in A, n\in \N\right\}\cr}$
\m
\item{iii)} \ \ $\pconv (T(A))= T(\pconv (A))$,
 for any linear map $T$.

\Pf i) and iii) are straighforward.

ii) We only have to prove that
$\pconv (A\cup\{ 0\})\subseteq \{ 0\}\cup \pconv (A)$. It is enough to show that
every non zero element  $x$ 
of the form $\d x=\lx,\ x_i\in A,\ \l< 1$ can be written as $\d x=\my$, $y_i\in
A$, 
 $\d \sum_1^m
\mu_i^p=1$. 

 Suppose $\lambda_1\ne 0$. Write $\d\lambda_1=\sk \beta_i$, with $
\beta_i\geq 0$. We have
$$\l \leq \sk\beta_i^p+\sum_{i=2}^n\lambda_i^p\leq
k^{1-p}\lambda_1^p+\sum_{i=2}^n\lambda_i^p.$$
 It is now clear, by a continuity
argument,  that we can find $k$ and $\beta_i\geq 0$, $1\le i \le k$, such that
$\d\lambda_1=\sk \beta_i$ and  $\d \sk \beta_i^p+\sum_{i=2}^n\lambda_i^p=1$.
Finally, the representation $\d x=\sk \beta_i x_i+\sum_{i=2}^n\lambda_i x_i$
 does the job. 
\fin

\Rm In particular ii) says that for every $0\ne x\in X $, $\pconv \{ x\}
=(0,x] = \{ \lambda x; 0<\lambda \leq 1\}$. This situation is rather different
from the case when $p=1$.

Another useful particular case of ii) is the following: If $\d A=\{P_1,\dots
,P_n\}\subset X, P_i\ne 0, P_i\ne P_j, \forall\,1\leq i\ne j\leq n$, then 
$\d 0\ne x\in\pconv(A)\Rightarrow x=\sum_{i=1}^n\lambda_i P_i$, with $ \l\leq
1, \lambda_i\geq 0$.  Observe that we allow  no more than $n$ non-zero
summands while in i) and ii) there is no restriction.

Next, we are going to prove a particular case of Theorem 1, which
will help us in the general case.

\proclaim Lemma 4. Let $\{ P_1 \ldots P_n, Q\} \subseteq \R^n$ with 
$\{P_i \}_{i=1}^n$ linealy independent. If $M \in \pconv {\{ P_1 \ldots P_n,
 Q, 0\}}$ then there exist $P_{i_1} \ldots P_{i_n} \in \{ P_1 \ldots P_n,
 Q\}$ such that $M \in \pconv{\{ P_{i_1} \ldots P_{i_n}, 0\}}$.

\Pf By Proposition 3 iii), it's enough to consider the case
$ \{ e_1 \ldots e_n, Q, 0\}$ where $\{e_i \}_{i=1}^n$ is the canonical basis in
$\R^n$ and $Q=(a_1 \ldots a_n) \not=0$. Denote by ${\cal P}$ the subset of
$\pconv\{ e_1 \ldots e_n, Q, 0\}$ for which the thesis of the Lemma holds.

Let $$K=\{ (\lambda_1 \ldots \lambda_n) \in \R^n\,\mid\,\l \le 1,\,\,\lambda_i
\ge 0, \,\,1\le i \le n\}$$
Write $\d \mu=\mu(\lambda_1 \ldots \lambda_n)=(1-\l)^{1/p},$
and consider the map $\varphi\colon K \to \R^n$ defined as $\d\lambda=
(\lambda_1 \ldots \lambda_n) \rightarrow \varphi(\lambda)={\sum_{i=1}
^n\lambda_i e_i}+\mu Q$. Denote by $\d {\rm J}\varphi(\lambda)$ the Jacobian
of the function $\varphi$ at a point $\lambda$. 

The proof of the Lemma rests on the following:

\proclaim Claim.  For every $\lambda\in{\rm Int}(K)$ such that $\d {\rm
J}\varphi(\lambda)=0$ we have $\varphi(\lambda)\in{\cal P}$.

Assume the claim is true and continue with the proof of the Lemma.

Let $ M \in \pconv{\{ e_1 \ldots e_n, Q, 0\}}$ i.e. $M={\sum_{i=1}^n
\lambda_i e_i}+\nu Q $, with $\lambda_i, \nu\geq 0$, $ 1\le i \le n$, $ 
\l+\nu^p \le 1$.

Suppose first that $\l+\nu^p =1$, that is $M=\varphi(\lambda), \lambda
\in K$. If $\lambda \in \partial(K) $ then, either $\lambda_i$ or $\nu$ are 
equal to $0$ and clearly $\varphi(\lambda)\in{\cal P}$. If $\lambda \in {\rm
Int}(K)$, we also have two posibilities: $a)$  ${\rm J}\varphi(\lambda)=0$ and
the claim says that $M \in {\cal P}$ or $b)$  ${\rm J}\varphi(\lambda)\ne 0$. 
By the inverse function theorem we necessarily have that $M\in {\rm Int}
\varphi(K)$. Since $\varphi(K)$ is compact 
there exists $ t>1$ such that $tM\in \partial\varphi(K) $, and therefore 
$tM= \varphi({\lambda}') $ with either ${\lambda}' \in \partial(K) $ or
$\lambda'\in{\rm Int}(K)$ and ${\rm J}\varphi(\lambda')=0$. In any case we
deduce that $tM$ belongs to ${\cal P}$ and so does $M$.

If, on other hand, $\d\l+\nu^p =s^p<1$ the results easily follows by considering
$\d{M\over s}$ and applying the preceding case.
\fin

{\bf Proof of the Claim}: It is an easy exercise to show that the Jacobian
of $\varphi$ is 
$$\d {\rm J}\varphi(\lambda)= 1-\sum_{i=1}^n a_i
\left({\mu \over \lambda_i}\right)^{1-p}.$$
 
For every $\lambda \in {\rm Int}(K)$ we write $\lambda=tv, v=(v_1 \ldots v_n),
v_i>0$, $ 1\le i \le n$, $ 0<t<1$, $\sum_{i=1}^n v_i^p =1$.
 We have ${\rm J}\varphi(\lambda)=0$ if and only if
 $$\left({\mu \over \lambda_i}\right)^{1-p}
\sum_{i=1}^n {a_i \over v_i^{1-p}} =1 \qquad {\rm  and}\qquad \mu^p=1-t^p$$
Write  $\d
R=\left(\sum_{i=1}^n{a_i \over v_i^{1-p}}\right)^{1 \over
 1-p} >0$.
 It is easy to see that for every $v$, 
there is a unique $t \in (0,1)$
such that ${\rm J}\varphi(tv)={\rm J}\varphi(\lambda)=0$. 
Explicitly $\d\lambda=
{R \over (1+R^p)^{1/p}}v$. 
Therefore with this new notation, the points
 $\lambda$ with ${\rm J}\varphi(\lambda)=0$
are such that
$\d M=\sum_{i=1}^n{Rv_i+a_i \over {(1+R^p)}^{1/p}}e_i$ 
where $ v_i>0$, $1\le i \le n$,
 $\d \sum_{i=1}^n v_i^p =1$, $\d R^{1-p}=\sum_{i=1}^n{a_i \over v_i^{1-p}}>0$.

\noindent{\bf  Case 1}. If $Rv_i+a_i \ge 0$ for all $i$, then $M \in
\pconv{\{
 e_1 \ldots e_n, 0\}}$.  Indeed, let's show that 
$$\d \sum_{i=1}^n{(Rv_i+a_i)}^p <
1+R^p$$
 This is equivalent to 
$$\d \sum_{i=1}^n\left(v_i+{a_i \over R}\right)^p -
\sum_{i=1}^n v_i^p -{1\over R} \sum_{i=1}^n{a_i \over v_i^{1-p}} < 0$$
and to 
$$\d \sum_{i=1}^n v_i^p {\left(1+{a_i \over R v_i}\right)}^p- \sum_{i=1}^n
v_i^p  \left(1+{a_i\over R v_i}\right) < 0.$$
 But this is obvious by the elementary inequality:
$$ (1+x)^p \le 1+px,\qquad x \ge -1 \eqno(*) $$

\noindent{\bf  Case 2}. If there is some $i, 1\le i \le n$ such that $Rv_i+a_i <
0$, then $a_i<0$. We suppose without loss of generality that $\d{\rm min}
\{{a_i\over  v_i}\mid 1\le i \le n \}$ is achieved at $i=1$. We shall prove 
$M \in \pconv{\{e_2 \ldots e_n, Q, 0\}}$. Recall that $\d M=\sum_{i=1}^n
{Rv_i+a_i \over (1+R^p)^{1/p}}e_i$. We will show that $M$ can be represented
as $\d M=\sum_{i=2}^n {\alpha}_i e_i+ \beta Q$ with $\d \sum_{i=2}^n {\alpha}_i^p
+\beta^p<1$. Comparing the two representations, it is easy to see that

$\d {\alpha}_i=\left(v_i-{v_1 a_i\over a_1}\right){R\over {(1+R^p)}^{1/p}}$,
$2 \le i\le n$

$\d \beta=\left(1+{Rv_1\over a_1}\right){1\over {(1+R^p)}^{1/p}}$.

By hypothesis we have $\beta$, ${\alpha}_i \ge 0$,  $2 \le i\le n$. It remains to
show that $\d \sum_{i=2}^n {\alpha}_i^p+\beta^p<1$, which is the same as
$$\d \sum_{i=2}^n {\left(v_i-{v_1 a_i\over a_1}\right)}^p+{\left(1+{Rv_1\over
a_1}\right)}^p{1\over R^p}< 1+{1\over R} \sum_{i=1}^n{a_i \over v_i^{1-p}}$$ or 
$$\d \left(\sum_{i=1}^n{a_i \over R v_i^{1-p}}\right){\left(1+{Rv_1\over
a_1}\right)}^p+ \sum_{i=2}^n v_i^p{\left(1-{v_1 a_i\over a_1
v_i}\right)}^p-\sum_{i=1}^n v_i^p\left(1+ {a_i \over Rv_i}\right)<0$$

 Again (*)  establishes
 $\d {\left(1+{Rv_1\over a_1}\right)}^p \le 1+p{Rv_1\over a_1}$ and
$\d {\left(1-{v_1 a_i\over a_1 v_i}\right)}^p \le 1-p{v_1 a_i\over a_1 v_i}$ and
the result easily follows.
\fin

We are now ready to state and prove the main theorem of the section.

\proclaim Theorem 1. Let $ A\subseteq \R^n$ and $0<p<1$. For every  $x\in \pconv
({A}), x\not=0$ there exist linearly independent vectors $\{ P_1 \ldots
P_k\}\subseteq A$  with $k\le n$, such that $x\in \pconv{\{ P_1 \ldots P_k\}}$. 
Moreover, if  $0\in \pconv({A}) $, there exits $\{ P_1 \ldots P_k\} \subseteq A$
with $k\le n+1 $ such that $0\in \pconv{\{ P_1 \ldots P_k\}}$.

\Pf Let  $x\in \pconv ({A})$, $x\not=0$, then $x=\sum_{i=1}^N {\lambda}_i P_i$
with $P_i \in A$, $ P_i \not=0$, $ \sum_{i=1}^N {\lambda}_i^p \le 1$,
${\lambda}_i>0$ and $1\le i \le N $. Let dim$( {\rm span}
\{P_i\}_{i=1}^N)=m\le n$. By Proposition 3 iii) and without loss of generality,
we can suppose that we are in $R^m$ and that $\d x=\sum_{i=1}^N {\lambda}_i P_i$
with $P_1 \ldots P_m$ linearly independent.

 Write $\d s^p=\sum_{i=1}^{m+1} {\lambda}_i^p$ and $\d \tilde x =
\sum_{i=1}^{m+1} {{\lambda}_i \over s} P_i$. Clearly $\tilde x \in \pconv {\{P_1
\ldots P_{m+1}\}}$ and therefore,  by Lemma 4 there exists $\{P_{k_1} \ldots
P_{k_m}\} \subset \{P_1 \ldots P_{m+1}\}$ such that $\d \tilde x =\sum_{i=1}^m
{\beta}_i \-P_{k_i}$, $\d\sum_{i=1}^m {\beta}_i^p \le 1$. Hence $$ x =
s\tilde x+\sum_{i=m+2}^N {\lambda}_i P_i=\sum_{i=1}^m s\beta_iP_{k_i}+
\sum_{i=m+2}^N {\lambda}_i P_i$$ with 
 $\d \sum_{i=1}^m {\beta}_i^p s^p + \sum_{i=m+2}^N {\lambda}_i^p
\le s^p+ \sum_{i=m+2}^N {\lambda}_i^p \le 1$.

We have represented $x$ as a combination of points of $A$ of length $N-1$.
Consider now, ${\rm span }\{\-P_{k_1} \ldots \-P_{k_m}, P_{m+2} \ldots P_N\}$
and repeat the argument until reaching a representation of length $\le n$.

If  $0\in \pconv({A}) $ then $\d 0=\sum_{i=1}^N\lambda_i P_i$,  $P_i\in A$,
${\lambda}_i>0$, $1 \le i \le N $ and $\d \sum_{i=1}^N {\lambda}_i^p =1$. As
before, we can suppose $P_1 \ldots P_m$ linearly independent with $m \le n$.
We consider $\d\sum_{i=1}^{m+1} {\lambda}_i P_i = -
\sum_{i=m+2}^N {\lambda}_i P_i  $. 
If we apply  Lemma 4 to $\d \tilde x=\sum_{i=1}^{m+1} {{\lambda}_i \over s} P_i
$, $\d s^p=\sum_{i=1}^{m+1} {\lambda}_i^p $ we obtain
$$\sum_{i=1}^m {\beta}_i P_i = - \sum_{i=m+2}^N {\lambda}_i P_i  $$
with $ \sum_{i=1}^m {\beta}_i^p \le 1$.
Hence  $0 \in p$-convex envelope of $N-1$ points.
Repeat the argument until reaching a representation of length $\le n+1$.
\fin
\vskip 12pt

\noindent{\bf 3. Gluskin's theorem for $0<p<1$.}
\vskip 12pt

In this section we are going to prove Theorem 2. As quoted above, Peck
showed that diam $\d ({\cal M}_n^p)\leq n^{2/p-1}$. Given an $n$-dimensional
$p$-normed space $X$, he considered its Banach envelope $X^b$ (the normed space
whose unit ball is the convex envelope of the unit ball of $X$) and proved
$d(X,X^b)\leq n^{1/p-1}$  (see [Pe] or [G-K]). By using John's theorem he
obtained the estimate. We want to prove that this result is sharp. More
precisely what we are going to show is \b 

\proclaim Theorem 2. Let  $0<p<1$. There exits a constant $C_p >0$ such that
for every $n\in \N $ 
$$C_p n^{2/p-1}\leq {\rm diam}({\cal M}_n^p)\leq n^{2/p-1}.$$

The proof of Theorem 2 follows Gluskin's original ideas. We first introduce
some notation. $S^{n-1}$ will denote the euclidean sphere in $\R^n$ with its
normalized Haar measure $\mu_{n-1}$ and $\Omega$ will be the product space 
$\d S^{n-1}\times\buildrel n)\over\dots\times S^{n-1}$ endowed with the product
probability $\P$.  If $K\subseteq\R^n$, $|K|$ is the Lebesgue measure of $K$. If
$A=(P_1,\dots ,P_n)\subset\Omega$, we write $\d Q_p(A)=\pconv\{\pm e_i, \pm
P_i\mid 1\leq i\leq n\}$, being $\{e_i\}_{i=1}^n$ the canonical basis of $\R^n$.
We denote by $\d \Vert\cdot\Vert_{Q_p(A)}$ the $p$-norm in $\R^n$ whose unit
ball is $ Q_p(A)$. 
\m

We only need to prove that for some absolute constant $C_p>0$, there exist
$A,A'\in\Omega$ such that simultaneously
$$\Vert T\Vert_{Q_p(A)\to Q_p(A')}\geq C_p n^{1/p-1/2}\qquad {\rm and}\qquad 
\Vert T^{-1}\Vert_{Q_p(A')\to Q_p(A)}\geq C_p n^{1/p-1/2}$$
hold for any $T\in {\rm SL}(n)$ (that is, any linear isomorphism
in $\R^n$ with det $T=1$). Straightforward argument shows that it is enough to
see that for any fixed $A'\in\Omega$ we have, 
$$\P\{\,A\in\Omega\mid \,\Vert T\Vert_{Q_p(A)\to Q_p(A')}<C_pn^{1/p-1/2}\, {\rm
for\ some\ } T\in {\rm SL}(n)\,\}<{1\over 2}$$

Fix $A\in\Omega$ and $t>0$. Consider the set 
$$\Omega(A',t)=\{\,A\in\Omega\mid \,\Vert T\Vert_{Q_p(A)\to Q_p(A')}<t\,
\,\,{\rm for\ some\ } T\in {\rm SL}(n)\,\}$$

The proof of the following three lemmas are analogous to the ones in the case
$p=1$ (see [Tmcz], \S 38).
\m

\proclaim Lemma 5. Let $A'\in\Omega$ and $t>0$. There exists a $t$-net $\d
N(A',t)$ in $\d \{\,T\in {\rm SL}(n)\, \mid \Vert T\Vert_{\ell_p^n\to
Q_p(A')}\leq t\}$ with respect to the metric induced by 
$\d \Vert\cdot\Vert^p_{\ell_2^n\to Q_p(A')}$ of cardinality 
$$|\,N(A',t)\,|\leq (3n^{1/p-1/2})^{n^2}{|Q_p(A')|^n\over |\,\{T\in 
{\rm SL}(n)\, \mid \Vert T\Vert_{\ell_2^n\to
\ell_2^n}\leq 1\}\,|}$$

\m

\proclaim Lemma 6. For every $A'\in\Omega$ and $t>0$ we have,
$$\Omega(A',t)\subseteq\bigcup_{T\in N(A',t)} \{\,A\in\Omega\mid \Vert
T(P_i)\Vert_{Q_p(A')}\leq 2^{1/p}t, \forall\,1\leq i\leq n\,\}$$

\m

\proclaim Lemma 7. Given $T\in  {\rm SL}(n)$, $A'\in\Omega$ and $t>0$,
$$\P\{\,A\in\Omega\mid \Vert T(P_i)\Vert_{Q_p(A')}\leq 2^{1/p}t, \forall\,1\leq i\leq
n\,\}\leq (2^{1/p}t)^{n^2}\left({|\ Q_p(A')\ |\over |B_{\ell_2^n}|}\right)^n$$

{\bf Proof of Theorem 2}: Numerical constants are always denoted by the same
letters $C$ (or $C_p$, if it depends only on $p$) although they may have
different value from line to line. Using consecutively the three preceding
lemmas we have for every $A'\in\Omega$ and $t>0$, $$\P\big(\Omega(A',t)\big)\leq
(C_ptn^{1/p-1/2})^{n^2}{|Q_p(A')|^{2n}\over |B_{\ell_2^n}|^n\cdot |\,\{T\in 
{\rm SL}(n)\, \mid \Vert T\Vert_{\ell_2^n\to \ell_2^n}\leq 1\}\,|}$$

It is well known that for some absolute constant $C>0$, (see [Tmcz]), 
$$|\,\{T\in  {\rm SL}(n)\, \mid \Vert T\Vert_{\ell_2^n\to
\ell_2^n}\leq 1\}\,|\geq C^{n^2}|B_{\ell_2^n}|^n$$

Now using Theorem 1 it is clear that if $A'=\{P_1,\dots P_n\}$, then
$\d Q_p(A')\subseteq\bigcup\pconv \{P_{k_1},\dots ,P_{k_n}\}$ where the union
runs over the $\d 4n\choose n$ choices of $\{P_{k_i}\}_{i=1}^n\subseteq\{\pm
e_i, \pm P_i, 1\leq i\leq n\}$. Since $\Vert P_i\Vert_2=1$ and 
$$ |\pconv\{P_{k_1},\dots ,P_{k_n}\}|= |{\rm det}\ [P_{k_1},\dots ,
P_{k_n}]\,|\cdot |\pconv\{e_1,\dots ,e_n\}|$$ we get
$$|Q_p(A')|\leq {4n\choose n}{|B_{\ell_p^n}|\over 2^n}\leq C_p^n n^{-n/p}2^{-n}$$
for some constant $C_p$ (see [Pi], pg 11). Hence, 
$$\P\big(\Omega(A',t)\big)\leq (C_ptn^{1/2-1/p})^{n^2}$$
If we take  a suitable $t>0$, we can assure $\d \P\big(\Omega(A',t)\big)<{1\over
2}$ and the result follows. 
\fin

\Rm In the same way as quoted above, given a $p$-normed space $X$ and $p<q\leq
1$, we can define the $q$-Banach envelope of $X$ as the $q$-normed space, $X^q$
whose unit ball es the $q$-convex envelope of the unit ball of $X$. It is easy
to see that $\d d(X,X^q)\leq d(X,Y)$ for any $n$-dimensional
$q$-normed space $Y$. Theorem 1 shows that $\d d(X,X^q)\leq
n^{1/p-1/q}$. Indeed, for every $\d x\in B_{X^q}=\qconv (B_X)$ and $\d \Vert
x\Vert_{X^q}=1$ there exist $P_1,\dots ,P_n\in B_X$ such that $x=
\s\lambda_i P_i$ with $\lambda_i\geq 0, 1\leq i\leq n, \s\lambda_i^q\leq 1$ and 
$$1\leq \Vert x\Vert_X\leq \s\lambda_i^p\Vert P_i\Vert_X^p\leq \l\leq
n^{1/p-1/q}$$
by homogeneity we achieve the result. Now it is easy to see that if $X,Y$ are
the spaces appearing in Theorem 2, then $\d d(X,X^q)\geq C_p n^{1/p-1/q}, 
\d d(Y,Y^q)\geq C_p n^{1/p-1/q}$ and $\d  d(X^q,Y^q)\geq C_p
n^{2/q-1}$.
In particular, for $q=1$,  $\d d(X,X^b)\geq C_pn^{1/p-1}$
, $\d d(Y,Y^b)\geq C_pn^{1/p-1}$ and $\d d(X^b,Y^b)\geq C_pn$.

\m

 \noindent {\bf Acknowledgments.} The authors are indebted to Yves Raynaud for
some comments in the proof of Lemma 4.
\b
\b
\centerline{\bf References.}

\b

\item{[Eg]} Eggleston, H.G.: Convexity. Cambridge Tracts in Math. and Math.
Phys.{\bf 47}. Cambridge University Press (1969). \m

\item{[Gl]} Gluskin, E.D.: The
diameter of the Minkowski compactum is approximately equal to $n$. Functional
Anal. and Appl. {\bf 15}(1), 72-73 (1981). \m

\item{[G-K]}  Gordon, Y., Kalton, N.J.: Local structure for quasi-normed spaces.
Preprint, (1992).
\m 

\item{[Pe]} Peck, T.: Banach-Mazur distances and projections on $p$-convex
spaces. Math. Zeits. {\bf 177}, 132-141 (1981). \m

\item{[Pi]} Pisier G.: The volume of convex bodies and Banach Space Geometry.
Cambridge University Press (1989).
 \m

\item{[Sz]} Szarek, S.J.: Volume estimates and nearly Euclidean decompositions
of normed spaces. S\'eminaire d'Analyse Fonctionnelle  \'Ecole Poly. Paris.
Expos\'e 25. (1979-80)\m

\item{[Tmzc]} Tomczak-Jaegerman, N.: Banach-Mazur distances and
finite-dimensional operator ideals. Pitman Monographs {\bf 38} (1989).

  \bye